\documentclass[12pt]{amsart}
\usepackage{url, calc, bm}
\usepackage{amsmath,amsxtra,amssymb,latexsym,epsfig,amscd,amsthm,fancybox,epsfig}
\usepackage[mathscr]{eucal}
\usepackage{graphicx}
\usepackage{multicol,xcolor}
\usepackage{epsfig} % for postscript graphics files
\usepackage{epstopdf}
\usepackage{cases}
\usepackage{subfig}
\usepackage{color}
\usepackage{hyperref}
\usepackage{bigints}
\usepackage{enumerate}
\usepackage{pdfpages}

\numberwithin{equation}{section}

\usepackage[left=3cm,right=3cm,top=3cm,bottom=3cm]{geometry}
\usepackage{amsfonts}
\usepackage{t1enc}
\usepackage{amssymb}
\usepackage{epsfig}
\usepackage[french,english]{babel}

\usepackage{amsmath}
\usepackage{graphics}
\usepackage{layout}
\usepackage{latexsym}
\usepackage{color}
\DeclareMathOperator{\Tr}{Tr}
\def\mylabel#1{\label{#1}}
\newtheorem{theorem}{Theorem}[section]
\newtheorem{lemma}[theorem]{Lemma}
\newtheorem{corollary}[theorem]{Corollary}
\newtheorem{proposition}[theorem]{Proposition}
\newtheorem{example}[theorem]{Example}
\newtheorem{remark}[theorem]{Remark}
\newtheorem{conjecture}[theorem]{Conjecture}

\newtheorem{definition}[theorem]{Definition}
\newtheorem{exercice}[theorem]{Exercice}

\def\bit{\begin{itemize}}

\def\eit{\end{itemize}}

\reversemarginpar   %%puts label in LEFT margin

\def\bc{\begin{center}}

\def\ec{\end{center}}

\def\bthm{\begin{theorem}}

\def\ethm{\end{theorem}}

\def\bcor{\begin{corollary}}

\def\ecor{\end{corollary}}

\def\bprop{\begin{proposition}}

\def\eprop{\end{proposition}}

\def\blem{\begin{lemma}}

\def\elem{\end{lemma}}

\def\bex{\begin{example}}

\def\eex{\end{example}}

\def\bexo{\begin{exercice} \rm }

\def\eexo{\end{exercice} }

\def\brem{\begin{remark}}

\def\erem{\end{remark}$\Box$}

\def\bdes{\begin{description}}

\def\edes{\end{description}}

\def\beq{\begin{equation}}

\def\eeq{\end{equation}}

\def\ben{\begin{enumerate}}

\def\een{\end{enumerate}}

\def\beqar{\begin{eqnarray}}

\def\eeqar{\end{eqnarray}}

\def\beqarr{\begin{eqnarray*}}

\def\eeqarr{\end{eqnarray*}}

\def\ZZ{{\mathbb Z}}       %bold Z
\def\RR{{\mathbb R}}  % bold R
\def\CC{{\mathbb C}}

\def\EE{{\mathbb E}}
\def\PP{{\mathbb P}}

\def\bx{{\mathbf x}}

\def\rar{\rightarrow}
\def\tX{\tilde X}

\def\eps{\varepsilon}

\def\1{{\rm 1\mskip-4.4mu l}}
\begin{document}
\title[Predator-prey with random switching]{On a predator--prey system with random switching that never converges to its equilibrium}
\author[A. Hening]{Alexandru Hening }
\thanks{A. Hening is in part supported by the NSF through the grant DMS 1853463.}
\address{Department of Mathematics\\
Tufts University\\
Bromfield-Pearson Hall\\
503 Boston Avenue\\
Medford, MA 02155\\
United States
}
\email{Alexandru.Hening@tufts.edu}
\author[E. Strickler]{Edouard Strickler}
\thanks{E. Strickler was in part supported by SNF grant 2000020-149871/1}
\address{ Institut de Math\'ematiques\\Universit\'e de Neuch\^atel, Switzerland}
\email{edouard.strickler@unine.ch}
%\bibliographystyle{amsplain}
%\bibliographystyle{apalike}
%\date{ }

\begin{abstract}
We study the dynamics of a predator-prey system in a random environment. The dynamics evolves according to a deterministic Lotka--Volterra system for an exponential random time after which it switches to a different deterministic Lotka--Volterra system. This switching procedure is then repeated. The resulting process is a Piecewise Deterministic Markov Process (PDMP). In the case when the equilibrium points of the two deterministic Lotka--Volterra systems coincide we show that almost surely the trajectory does not converge to the common deterministic equilibrium. Instead, with probability one, the densities of the prey and the predator oscillate between $0$ and $\infty$. This proves a conjecture of Takeuchi et al (J. Math. Anal. Appl 2006).

The proof of the conjecture is a corollary of a result we prove about linear switched systems. Assume $(Y_t, I_t)$ is a PDMP that evolves according to $\frac{dY_t}{dt}=A_{I_t} Y_t$
where $A_0,A_1$ are $2\times2$ matrices and $I_t$ is a Markov chain on $\{0,1\}$ with transition rates $k_0,k_1>0$. If the matrices $A_0$ and $A_1$ are not proportional and are of the form
\[
A_i := \begin{pmatrix}
   \alpha_i & \beta_i \\
   \gamma_i & -\alpha_i
\end{pmatrix},
\]
with $\alpha_i^2 + \beta_i \gamma_i < 0$, then there exists $\lambda >0$ such that $\lim_{t \to \infty} \frac{\log \| Y_t \|}{t} = \lambda.$
\end{abstract}

\maketitle

%\newpage
%\tableofcontents
%\newpage
\section{Introduction and main results}
\mylabel{sec:intro}
One of the key issues in ecology is determining when species will persist and when they will go extinct.  The randomness of the environment makes the dynamics of populations inherently stochastic and therefore we need to take into account the combined effects of biotic interactions and environmental fluctuations. One way of doing this is by modelling the densities of various species as Markov processes and looking at their long-term behavior (see \cite{C00, ERSS13, EHS15,  LES03, SLS09, SBA11, BEM07, BS09, BHS08, CM10, HNY17, HN16}).

In order to allow for environmental fluctuations and their effect on the persistence or extinction of species one approach is to study stochastic differential equations (\cite{ERSS13, SBA11, HNY17, HN16, HN17, HN17b}). The other possible approach is to look at stochastic equations driven by a Markov chain. These systems are sometimes called\textit{ Piecewise Deterministic Markov Processes (PDMP)} or systems with \textit{telegraph noise}.

PDMPs have been used recently to prove some very interesting and counterintuitive facts about biological populations. In \cite{BL16} the authors look at a two dimensional competitive Lotka--Volterra system in a fluctuating environment. They show that the random switching between two environments that are both favorable to the same species can lead to the extinction of this favored species or to the coexistence of the two competing species (also see \cite{MalHoa16}). PDMPs are also used in \cite{M16} where the author studies prey-predator communities where the predator population evolves much faster than the prey.

For a predator-prey system the classical deterministic example is the Lotka--Volterra model (see \cite{L25} and \cite{V28})
\begin{equation}\label{e:LV}
\begin{split}
\frac{dx(t)}{dt} &= x(t)(a-by(t)),\\
\frac{dy(t)}{dt} &= y(t)(-c+dx(t)),
\end{split}
\end{equation}
where $x(t), y(t)$ are the densities of the prey and the predator at time $t\geq 0$ and $a,b,c$ and $d$ are positive constants. If one assumes that $x(0)=x_0>0, y(0)=y_0>0$, so that both predator and prey are present, then the solutions of system \eqref{e:LV} are periodic (see \cite{G75, HS98}) and given in phase space by the curves described by the first integral,
\begin{equation}\label{e:first}
r(x,y)=dx-c-c\ln(1+(dx-c)/c) + by-a-a\ln(1+(by-a)/a)=\text{constant}=r.
\end{equation}
One should note that both the predator and the prey from \eqref{e:LV} do not experience intraspecific competition. In particular, if the predator is not present (i.e. $y_0=0$) then the prey density blows up to infinity. In \cite{GH79, MHP14} the authors are able to analyze the $n$-dimensional generalization of \eqref{e:LV} i.e. the setting when one has one prey and $n-1$ predators and each species interacts only with the adjacent trophic levels. Stochastic predator-prey models have been studied in the stochastic differential equation setting by \cite{R03, HN17, HN17b}. However, we note that in all these studies one needed to assume that there exists intraspecific competition among the prey and the predators. This simplifies the analysis significantly because the predator and the prey densities get pushed towards the origin when they become too large.

In \cite{AHS79} the authors show that if the coefficient $a$ (growth rate of the prey) is randomly perturbed by white noise then the resulting stochastic system cannot have a stationary distribution and that as the time goes to infinity, with probability $1$, explosion does not occur. In \cite{KK01} the authors look at scaling limits of Lotka--Volterra systems perturbed by white noise - they prove that a suitably rescaled version of $r(x(t),y(t))$, where $r(x,y)$ is the first integral from \eqref{e:first}, converges to a one-dimensional stochastic differential equation. They then use this SDE to gain information about both the deterministic and the stochastic Lotka--Volterra systems.

We consider the random switching between two Lotka--Volterra prey-predator systems of the form \eqref{e:LV}. More precisely, for $i \in E:=\{0,1\}$, let $F^i:\RR_+^2\to\RR_+^2$ denote the vector field
\begin{equation}\label{e:vector_field}
F^i(x,y) =  \begin{pmatrix}
x(a_i-b_iy)\\y(-c_i+d_ix)
\end{pmatrix}
\end{equation}
with $a_i, b_i, c_i, d_i > 0$.  Let $(I_t)_{t \geq 0}$ be a continuous-time Markov Chain defined on some probability space $(\Omega, \mathcal{F}, \PP)$ and taking values in $E:=\{0,1\}$. Suppose $I_t$ has transition rates $k_0, k_1 >0$.
Throughout the paper we will let $\RR_{++}^2:=\{(x_1,x_2)\in \RR^2 ~|~ x_1>0,x_2>0\}$ and  $\RR_{+}^2:=\{(x_1,x_2)\in \RR^2 ~|~ x_1\geq 0,x_2\geq 0\}$.
We denote by $(X_t)_{t \geq 0}=(x_t,y_t)_{t \geq 0}$ the solution of
\begin{equation}\label{e:LV_PDMP}
\begin{split}
\frac{dx_t}{dt}&= x_t(a_{I_t}-b_{I_t}y_t)\\
\frac{dy_t}{dt}&= y_t(-c_{I_t}+d_{I_t}x_t)
\end{split}
\end{equation}
for some initial condition $X_0=(x_0,y_0) \in \RR_{++}^2$. The process $(X,I)=(X_t, I_t)_{t \geq 0}$ is a Piecewise Deterministic Markov Process  as introduced in \cite{Dav84}, and belongs to the more specific class of PDMPs recently studied in \cite{BH12} and \cite{BMZIHP}.

One can construct the process $(X, I)$ as follows: Suppose we start at $(X_0, I_0)=((x_0,y_0), i)$. Then, the system evolves according to
\begin{equation*}
\begin{split}
\frac{dx_i(t)}{dt} &= x_i(t)(a_i-b_iy_i(t)),\\
\frac{dy_i(t)}{dt} &= y_i(t)(-c_i+d_ix_i(t))\\
x_i(0)&= x_0\\
y_i(0) &= y_0
\end{split}
\end{equation*}
for an exponential random time $T_i$ with rate $k_i$. After this time the Markov chain $I$ jumps from state $i$ to state $j\in\{0,1\}\setminus \{i\}$ and $X_t$ evolves according to
\begin{equation*}
\begin{split}
\frac{dx_j(t)}{dt} &= x_j(t)(a_j-b_jy_j(t)),\\
\frac{dy_j(t)}{dt} &= y_j(t)(-c_j+d_jx_j(t))\\
x_j(T_i)&= x_i(T_i)\\
y_j(T_i) &= y_i(T_i)
\end{split}
\end{equation*}
for an exponential random time $T_j$ with rate $k_j$. This procedure then gets repeated. Intuitively our process follows an ODE for an exponential random time after which it switches to a different ODE, follows that one for an exponential random time and so on.

The generator $L$ of $(X, I)$ acts on  functions $g : \RR_+^2 \times E \to \RR$ that are smooth in the first variable as $$Lg(x,i) = \langle F^i(x), \nabla g^i(x) \rangle + k_i \left( g(x,1-i) - g(x,i) \right),$$
where $\langle \cdot, \cdot \rangle$ is the euclidean inner product on $\RR^2$. As usual, for $\bx \in \RR^2$ and $i \in E$, we denote by $\PP_{\bx,i}$ the law of the process $(X,I)$ when $(X_0, I_0) = (\bx, i)$ almost surely and by $\EE_{\bx,i}$ the associated expectation.

The vector field $F^i$ from \eqref{e:vector_field} has a unique positive equilibrium $(p_i,q_i) = (c_i/d_i, a_i/b_i)$. In \cite{t06} the authors look at the two cases
\begin{itemize}
\item[Case I.] $p_0=p_1=:p$ and $q_0 = q_1=:q$, i.e. common zero for $F^0$ and $F^1$,
\item [Case II.] $(p_0,q_0)\neq (p_1,q_1)$, i.e. different zeroes for $F^0$ and $F^1$.
\end{itemize}

We assume throughout this paper that $p_0=p_1=:p$ and $q_0 = q_1=:q$. The vector fields $F^0$ and $F^1$ therefore have a common zero - this will allow us to use the recent results from \cite{BS17}. We also assume that $F^0$ and $F^1$ are non collinear to avoid trivial switching.

In \cite[Theorem 4.5]{t06} it is shown that only two long term behaviours are possible when the vector fields have a common zero: either $X_t$ converges almost surely to the common equilibrium $(p,q)$, or each coordinate oscillates between $0$ and $+ \infty$.

\bthm[Takeuchi et al., 2006]
For any $(x_0,y_0) \in \RR_{++}^2 $, with probability $1$, either
\beq
\lim_{ t \to \infty} X_t = (p,q),
\label{cveq}
\eeq
or
\beq
\limsup x_t = \limsup y_t = +\infty, \quad \liminf x_t = \liminf y_t =0.
\label{div}
\eeq
\label{thm:taku}
\ethm

It was conjectured from simulations (see \cite[Remark 5.1]{t06}) that only case \ref{div} happens in the above theorem. Using Theorem \ref{mainprop} below and results from \cite{BS17}, we are able to prove this conjecture.

\bthm\label{mainthm}
There exist $\eps > 0$,  $\eta > 1$, $\theta > 0$ and $C > 0$ such that
for all $\bx:=(x_0,y_0)  \in \RR_{++}^2\setminus \{(p,q)\}$ and $i \in E$, $$\mathbb{E}_{\bx,i}(\eta^{\tau^{\eps}}) \leq C (1 + \|\bx-(p,q)\|^{-\theta}),$$ where
\begin{equation}\label{e:tau}
\tau^{\eps} := \inf \{t \geq 0 : \: \|X_t-(p,q)\| \geq \eps \}.
\end{equation}
In particular, for any $(x_0,y_0)  \in \RR_{++}^2\setminus \{(p,q)\}$ we have with probability $1$ that
\[
\limsup_{t\to\infty} x_t = \limsup_{t\to\infty} y_t = +\infty, \quad \liminf_{t\to\infty} x_t = \liminf_{t\to\infty} y_t =0.
\]

\ethm

Our result provides a deeper understanding of Lotka--Volterra systems in random environments, continuing the work started in \cite{KK01} and \cite{AHS79}.

\subsection{Linear switched systems}\label{s:lin}
Let $\tilde A_i$ denote the Jacobian matrix of the vector field $F^i$ at $(p,q)$, where $(p,q)$ is the common positive equilibrium of the vector fields $F^0$ and $F^1$. Then $$\tilde A_i = \begin{pmatrix}
   0 & -b_i p \\
   d_i q & 0
\end{pmatrix} = \begin{pmatrix}
   0 & \beta_i \\
   \gamma_i  & 0
\end{pmatrix} ,$$
where $\beta_i = - b_i p$ and $ \gamma_i = d_i q$.

The matrices $\tilde A_i$ represent the linearizations of the nonlinear Lotka--Volterra systems near their common equilibrium point $(p,q)$. In order to study the dynamics of the nonlinear switched system \eqref{e:LV_PDMP} we will first study the linearization with switching and then use results of \cite{BS17}. Since we can prove slightly more general results for the linear systems we will work in the following setting.

Let $A_i$ denote the matrix
\begin{equation}\label{e:matrix}
A_i := \begin{pmatrix}
   \alpha_i & \beta_i \\
   \gamma_i & -\alpha_i
\end{pmatrix},
\end{equation}
for $i=0,1$, where $\alpha_i, \beta_i, \gamma_i$ are real numbers satisfying
\begin{equation}\label{e:imaginary}
\alpha_i^2 + \beta_i \gamma_i < 0.
\end{equation}
In this case, both matrices $A_0, A_1$ have purely imaginary eigenvalues.

% $$A_i := \begin{pmatrix}
%   0 & -\alpha_i \\
%   \beta_i & 0
%\end{pmatrix},$$
%for $i=0,1$, where $\alpha_i, \beta_i$ are positive.
We consider a random switching between the two dynamics given by $A_0$ and $A_1$. Let $(I_t)_{t \geq 0}$ be a continuous-time Markov Chain on $E=\{0,1\}$ with transition rates $k_0, k_1 >0$. We denote by $(Y_t)_{t \geq 0}$ the solution of
\begin{equation}\label{e:Y}
\begin{split}
\frac{dY_t}{dt}&=A_{I_t} Y_t\\
Y_0&=\mathbf{y}_0 \in \RR^2 \setminus \{(0,0)\}.
\end{split}
\end{equation}
 The process $(Y_t, I_t)_{t \geq 0}$ is a PDMP living on $\RR^2 \setminus \{(0,0)\} \times E$.

We will show that, independent of the starting conditions, $\|Y_t\|$ converges exponentially fast to infinity with probability one. More precisely, we prove the following.
\bthm\label{thm:growthrate}
Assume $A_0$ and $A_1$ are non proportional matrices of the form \eqref{e:matrix} with coefficients satisfying \eqref{e:imaginary}. Then, there exists $\lambda > 0$ such that, for all $\mathbf{y}_0 \neq 0$, almost surely
\beq
 \lim_{t \to \infty} \frac{1}{t} \log \| Y_t \| = \lambda.
 \label{e:lim}
\eeq

\ethm
\subsection{Generalization to density-dependent switching rates}
Actually, thanks to Theorem \ref{mainprop}, the first part of Theorem \ref{mainthm} can be significantly generalised.

For $i \in E$, let $F^i$ be a vector field of class $C^2$ on $\RR^2$, such that $F^i(0)=0$. Also assume that for $i \in E$, $DF^i(0)$, the Jacobian matrix of $F$ at 0, has two purely imaginary eigenvalues. In this case, the equilibrium 0 is sometimes called a \textit{center}. We now consider a Markov process $(U_t,J_t)_{t \geq 0}$ where $U_t$ is solution of $$\frac{dU_t}{dt}=F^{J_t} (U_t)$$ and $I_t$ is a jump process on $E$ whose rates depend on $X$ $$ \PP( J_{t+s} = 1 - i \vert J_t=i, \mathcal{F}_t) = k_{i, 1-i}(U_t)s + o(s),$$ where $F_t = \sigma( (U_s,J_s) : \: s \leq t)$ and for all $x$, $(k_{ij}(x))_{i,j}$ is an irreducible matrix that is continuous in $x$. The process $(U,J)$ is still a PDMP, with infinitesimal generator $\mathcal{L}$ acting  on  functions $g : \RR_+^2 \times E \to \RR$ that are smooth in the first variable as $$\mathcal{L}g(x,i) = \langle F^i(x), \nabla g^i(x) \rangle + k_{i, 1-i}(x) \left( g(x,1-i) - g(x,i) \right).$$
We can prove the following (see Remark \ref{r:extend}) in this more general setting.
\bthm \label{thm:gen}
Assume $DF^0(0)$ and $DF^1(0)$ are non proportional matrices such that, for $i \in \{0,1\}$,
$$
\Tr(DF^i(0)) = 0 \quad \mbox{and} \quad \det(DF^i(0)) > 0.
$$ Then there exist $\eps > 0$,  $\eta > 1$, $\theta > 0$ and $C > 0$ such that
for all $\bx:=(x_0,y_0)  \in \RR^2\setminus \{(0,0)\}$ and $i \in E$, $$\mathbb{E}_{\bx,i}(\eta^{\tau_U^{\eps}}) \leq C (1 + \|\bx\|^{-\theta}),$$ where  $\tau_U^{\eps} = \inf \{t \geq 0 : \: \|U_t\| \geq \eps \}.$
In particular, for any $(x_0,y_0)  \in \RR^2\setminus \{(0,0)\}$, with probability one, $U_t$ cannot converge to $(0,0)$.
\ethm

The paper is structured as follows. In Section \ref{sec:linear} we prove results about the linear switched systems introduced in Section \ref{s:lin}. In particular, we prove Theorem \ref{thm:growthrate}. We then apply these results in Section \ref{s:main} where we prove Theorems \ref{mainthm} and \ref{thm:gen}. Finally, in Section \ref{s:future} we present some conjectures and directions for future work.

\section{A result on linear switched systems}
\mylabel{sec:linear}
In this section we work with the linear systems introduced in Section \ref{s:lin} by equations \eqref{e:matrix}, \eqref{e:imaginary} and \eqref{e:Y}.
In order to do this we will use a \textit{polar decomposition}. The use of polar decompositions to study Lyapunov exponents goes back to \cite{Has60} in the case of stochastic differential equations. They have been used recently in the study of linear PDMPs (see \cite{BLMZexample,Lawley&matt&Reed}) and more general PDMPs (see \cite{BS17}).

Throughout the paper, we will denote by $S^1\subset \RR^2$ the circle with center at $0$ and radius $1$.
Whenever $\mathbf{y}_0 \neq 0$ and $Y_t\neq 0$, setting $\Theta_t= Y_t / \| Y_t\|$ and $\rho_t = \|Y_t\|$, one can check using \eqref{e:Y} that $(\rho_t,\Theta_t)_{t \geq 0}$ is the solution to
 \begin{equation}\label{e:polar}
\begin{split}
\frac{d\Theta_t}{dt}& = A_{I_t} \Theta_t - \langle A_{I_t} \Theta_t, \Theta_t \rangle \Theta_t\\
\frac{d\rho_t}{dt}& = \rho_t\langle A_{I_t} \Theta_t, \Theta_t \rangle \\
\Theta_0&=\theta_0 \in S^1\\
\rho_0 &= r_0 > 0,
\end{split}
\end{equation}
with $\theta_0 = \mathbf{y}_0/ \| \mathbf{y}_0\|$ and $r_0 = \| \mathbf{y}_0 \|$. In particular, $((\Theta_t,I_t))_{t \geq 0}$ is a PDMP on $S^1 \times E$ (see \cite{BS17}) and one has for all $t \geq 0$,
\begin{equation}
 \frac{1}{t} \log \| Y_t \| = \frac{1}{t} \int_0^t  \langle A_{I_s} \Theta_s, \Theta_s \rangle ds + \frac{1}{t} \log \| \mathbf{y}_0 \|.
 \label{eq:ergodicgrowth}
\end{equation}
Moreover, we have the following result.

\blem
Assume $A_0$ and $A_1$ are two matrices of the form \eqref{e:matrix} with coefficients satisfying \eqref{e:imaginary}. Then, the process $(\Theta_t,I_t)$ admits a unique invariant probability measure $\mu$ on $S^1 \times E$. Furthermore,
$$
\Lambda= \int \langle A_i \theta, \theta \rangle \mu(d \theta di) \geq 0.$$
\label{lem:growthnonnegative}
\elem

\begin{proof}
The uniqueness follows from \cite[Proposition 2.11 and Example 2.12]{BS17}. Indeed, since we study a two dimensional system, a sufficient condition is that at least one matrix $A_i$ has no real eigenvalue. This is the case for both $A_0$ and $A_1$. Since $A_0$ and $A_1$ have zero trace, \cite[Corollary 2.7]{BS17} implies that $\Lambda\geq 0$.
\end{proof}

By the Ergodic Theorem, \eqref{eq:ergodicgrowth} and Lemma \ref{lem:growthnonnegative}, one has for all $ \mathbf{y}_0 \neq 0$ and all $i \in E$
$$
\lim_{t \to \infty} \frac{1}{t} \log \| Y_t \| = \Lambda.
$$
Because of this, Theorem \ref{thm:growthrate} is a consequence of the following theorem, which is the main result of this section.

\bthm
Assume $A_0$ and $A_1$ are non proportional matrices of the form \eqref{e:matrix} with coefficients satisfying \eqref{e:imaginary}.
Then, with the notation of Lemma \ref{lem:growthnonnegative},
$$
\Lambda= \int \langle A_i \theta, \theta \rangle \mu(d \theta di) > 0.
$$
\label{mainprop}
\ethm

\subsection{Lyapunov exponents and Bougerol's theorem}

In order to prove Theorem \ref{mainprop}, we will use results from \cite{B88} on Lyapunov exponents. These numbers give the exponential growth rate of a linear \textit{Random Dynamical System} (see Arnold \cite{Arn98} for the definition). In \cite{BS17}, the authors show that the process $Y$   from \eqref{e:Y} together with the canonical shift on the space $\Omega$ of cadlag functions from $\RR_+$ to $E$ is a linear ergodic random dynamical system satisfying the integrability conditions of Osedelet's Multiplicative Ergodic Theorem (see \cite[Theorem 3.4.1]{Arn98} or \cite[Proposition 3.12]{ColMaz}).
According to this theorem, there exist $1\leq d\leq 2$ numbers  such that if $d=2$ $$ \lambda_1 > \lambda_2 $$ called {\em the Lyapunov exponents}, a Borel set $\tilde{\Omega}  \subset \Omega$ with $\mathbb{P}(\tilde{\Omega}) = 1,$ and for each $\omega \in \tilde{\Omega}$ distinct vector spaces
$$\{0\} = V_{d+1}(\omega) \subset V_{d}(\omega) \subset V_1(\omega) = \RR^2$$ such that
\beq
\label{eq:deflambdai}
\lim_{t \rar \infty} \frac{1}{t} \log \|Y_t\| = \lambda_i
\eeq for all $y_0 \in V_i(\omega) \setminus V_{i+1}(\omega)$. Now, since $\Tr(A_i)=0$ for $i \in E$, \cite[Corollary 2.7]{BS17} implies that
$$
\sum_{i} \lambda_i = 0.
$$
In addition, \cite[Proposition 2.5]{BS17} yields
$$
\Lambda = \lambda_1.
$$
Therefore, proving that $\Lambda > 0$ is equivalent to showing that $d=2$. In order to prove that $d=2$, we will use \cite[Theorem 1.7]{B88} which we state below (see Theorem \ref{thm:boug}).

Let $(M_t)_{t \geq 0}$ with $M_t\in  GL_2(\RR), t\geq 0$ be the solution of the matrix equation
 \begin{equation}\label{e:M_t}
 \begin{split}
 \frac{dM_t}{dt}&=A_{I_t} M_t\\
 M_0&\in  GL_2(\RR).
 \end{split}
 \end{equation}
Here $GL_2(\RR)$ stands for the set of invertible $2 \times 2$ matrices with real coefficients and the process $(M_t, I_t)$ is a PDMP living on $GL_2(\RR)\times E$.

The fact that $M_t \in GL_2(\RR)$ can easily be shown recursively as follows. If $T_1$ denotes the first jump time of the process $I$, then for all $t \leq T_1$, one has $M_t = e^{t A_{I_0}} M_0$, which is invertible as a product of invertible matrices. Then if $T_2$ is the second jump time, for all $t \in (T_1, T_2]$, $M_t = e^{(t-T_2) A_{1 - I_0}}M_{T_1}$, and so on.  One can note that when $M_0 = \mathrm{Id}$, the identity matrix, then for all $\mathbf{y} \in \RR_+^2$ the process $Y$ from \eqref{e:Y} can  be written as

$$Y_t = M_t \mathbf{y}$$
if $Y_0 = \mathbf{y}$.
Set $\pi = (k_1/(k_0+k_1), k_0/(k_0+k_1))$ and $\mathcal{E}=E=\{0,1\}$. Then $(M,I,\pi)$ is a simple example of a what Bougerol calls a \emph{multiplicative system} (see Definition \ref{def:multsys} and Lemma \ref{lem:multi} below).

We next present the abstract framework of \cite{B88}. Let $(\sigma_t)_{t \geq 0}$ be a stationary Markov process on some metric space $\mathcal{E}$, and $(B_t)_{t \geq 0}$ a process with values in $GL_2(\RR)$. We recall that the semigroup of a Markov process $(B,\sigma)$ is a family of measures defined for $t \geq 0$ and $(A,i) \in  GL_2(\RR) \times \mathcal{E} $ by $$P_t\left((A,i); \cdot\right) = \PP_{A,i} \left( (B_t,\sigma_t) \in \cdot \right).$$ Equivalently, the semigroup can be seen as a family $(P_t)_{t\geq 0}$ of operators which act on bounded measurable functions $f : GL_2(\RR) \times \mathcal{E} \to \RR$ according to $$P_t f((A,i)) = \EE_{A,i} \left[ f((B_t,\sigma_t)) \right], ~A\in GL_2(\RR), i\in \mathcal{E}.$$ We introduce the following definition, which is  \cite[Definition 1.1 and 1.2]{B88}.

\begin{definition}
 Let $\pi$ be a probability measure on $\mathcal{E}$. We say that $(B,\sigma,\pi)$ is a \emph{multiplicative system}, if the following properties hold
\begin{itemize}
\item [i)] The process $(B,\sigma)$ is Markovian with semigroup $(P_t)_{t \geq 0}$;
\item [ii)] For any Borel subset $A\subset \mathcal{E}$ (resp. $B\subset GL_2(\RR)$), $t \geq 0$, $C \in GL_2(\RR)$ and $i \in \mathcal{E}$, one has $$P_t\left((C,i);BC \times A \right) = P_t\left((\mathrm{Id},i);B \times A \right),$$
where $BC= \{ NC; N \in B \}$;
\item [iii)] $\pi$ is an ergodic measure for $\sigma$ and $ \sup_{0 \leq t \leq 1} \EE_{\mathrm{Id},\pi} \left[ \log^+ \| B_t \| + \log^+ \| B_t^{-1} \| \right] < \infty$.
\end{itemize}
\label{def:multsys}
\end{definition}

We let $Q$ be the first order resolvent of the semigroup $(R_t)_{t\geq 0}$ of the Markov process $\sigma$.  This is defined via $$ Q = \int_0^{+\infty} \mathrm{e}^{-t} R_t \mathrm{d} t.$$

\begin{remark}
The resolvent has the following properties related to the dynamics of its associated Markov process:
\begin{itemize}
  \item  A probability measure is invariant for the semigroup if and only if it is invariant for the resolvent.
\item  A point $y$ is accessible for the process from $x$ (meaning that it can reache every neighborhood of $y$ with strictly positive probability) if and only if $y$ is in the support of the resolvent
\end{itemize}
\end{remark}

Following  \cite{B88}, we will say that the semigroup $(P_t)_{t\geq 0}$ is \textit{Feller} if for any bounded continuous map $f : GL_2(\RR) \times \mathcal{E} \to \RR$, and for all $t \geq 0$, the function $P_t f$ is continuous.

\begin{definition}
We say that a multiplicative system $(B,\sigma,\pi)$ satisfies hypothesis \textbf{H} if the following conditions hold
\begin{itemize}
\item [i)] The space $\mathcal{E}$ is a complete metric space.
\item [ii)] The semigroup $(P_t)_{t \geq 0}$ is Feller.
\item [iii)] The support of $\pi$ is $\mathcal{E}$. If $h$ is a bounded measurable function which is a fixed point for the first order resolvent of $\sigma$, i.e.
\[
Qh= h
\]
then $h$ is continuous.
\label{def:H}
\end{itemize}
\end{definition}
Denote by $U$ the first order resolvent of $(P_t)_{t\geq 0}$, that is $$ U = \int_0^{+\infty} \mathrm{e}^{-t} P_t \mathrm{d} t.$$ For $i \in \mathcal{E}$ let $D_i$ be the support of $U((\mathrm{Id},i),\cdot)$ and $\mathcal{S}_i = \{ A \in GL_2(\RR) \: : \: (A,i) \in D_i \}$.
One has the following result (see \cite[Theorem 1.7]{B88}).
\bthm[Bougerol, 1988]\label{thm:boug}
Assume $(B,\sigma,\pi)$ is a multiplicative system defined on the space $\Omega$ of functions from $\RR_+$ to $GL_2(\RR) \times \mathcal{E}$ and satisfying hypothesis \textbf{H}. Assume furthermore that
\begin{itemize}
\item [i)] For some $i \in \mathcal{E}$, there exists a matrix in $\mathcal{S}_i$ with two eigenvalues with different modulus.
\item [ii)] There does not exist some finite union $W$ of one-dimensional vector spaces  such that, for all matrices $B$ in $\mathcal{S}_i$, $B W= W$.
\end{itemize}

Then there exist $2$ numbers  $ \lambda_1 > \lambda_2 $,  a Borel set $\tilde{\Omega}  \subset \Omega$ with $\mathbb{P}_{\mathrm{Id}, \pi}(\tilde{\Omega}) = 1,$ and for each $\omega \in \tilde{\Omega}$ distinct vector spaces
$$\{0\} = V_{3}(\omega) \subset V_{2}(\omega) \subset V_1(\omega) = \RR^2$$ such that
$$
\lim_{t \rar \infty} \frac{1}{t} \log \|B_t(\omega) y_0\| = \lambda_i
$$
 for all $y_0 \in V_i(\omega) \setminus V_{i+1}(\omega)$.
\ethm

\begin{remark}
Theorem \ref{thm:boug} is a reformulation of \cite[Theorem 1.7]{B88}, which is given for the numbers $\gamma_i$ that are the Lyapunov exponents for the external power of $M$ (see \cite[Proposition 2.2]{B88} or \cite[Theorem 3.3.3]{Arn98} for details). The numbers $\gamma_1$ and $\gamma_2$ are the numbers $\lambda_i$ counted with multiplicity. (see \cite[Definition 3.3.8 and Theorem 3.4.1]{Arn98}).
\end{remark}

 We show that we can use Theorem \ref{thm:boug} in our context.  Let $(M_t)_{t \geq 0}$ with $M_t\in  GL_2(\RR), t\geq 0$ be the process defined above, i.e. the solution of the matrix equation
 \begin{equation*}
 \begin{split}
 \frac{dM_t}{dt}&=A_{I_t} M_t\\
 M_0&\in  GL_2(\RR).
 \end{split}
 \end{equation*}

\blem
Set $\pi =  (k_1/(k_0+k_1), k_0/(k_0+k_1))$ and $\mathcal{E}=E=\{0,1\}$. Then $(M,I,\pi)$ is a multiplicative system satisfying \textbf{H}.
\label{lem:multi}
\elem

\begin{proof}
First we show that $(M,I,\pi)$ is a multiplicative system. $(M,I)$ is a PDMP thus a Markov process. In addition, if we denote by $M^N$ the process $M$ when $M_0 = N$ almost surely, one can easily check that $M^N = M^{\mathrm{Id}} N$ almost surely. As a result we see that point \textbf{(ii)} of definition \ref{def:multsys} is satisfied.
Straightforward computations show that $\pi$ is the unique invariant distribution of $I$ and is therefore ergodic.
 Let $K$ be a constant such that $\|A_i\| \leq K$ for $i \in E$. Then from $$\frac{d{M}_t}{dt} = A_{I_t} M_t,$$ $$\frac{d{M}^{-1}_t}{dt} = - M^{-1}_t A_{I_t}$$ together with $M_0=M_0^{-1}=\mathrm{Id}$ and Gronwall's Lemma, one can show that for all $t\geq 0$  $$\|M_t\|, \|M^{-1}_t\| \leq e^{K t}. $$This proves point \textbf{(iii)}.

Now we show that  $(M,I,\pi)$ satisfy hypothesis \textbf{H}. In our case, $\mathcal{E}=E$ is a finite set, thus points (i) and (iii) of Definition \ref{def:H} are straightforward. To prove that $(M,I)$ is Feller, we use \cite[Proposition 2.1]{BMZIHP} where the authors show that for a PDMP remaining in a compact set, the semigroup maps every continuous function to a continuous function.  Their proof adapts verbatim to the case where the process does not remain in a compact set with the additional assumption that the continuous function is bounded and provided the jump rates are bounded - in our setting the jump rates are constant. This concludes the proof.
\end{proof}

\subsection{Proof of Theorem \ref{mainprop}}

We start by showing that it suffices to prove Theorem \ref{mainprop} for a specific class of matrices.

\blem
Assume Theorem \ref{mainprop} holds when $A_0$ is of the special form $$ A_0 = \begin{pmatrix}
   0 & -\omega_0 \\
   \omega_0 & 0
\end{pmatrix}.$$
\label{lem:part}
Then Theorem \ref{mainprop} holds for any $A_0$.
\elem

\begin{proof}
First we show that a linear change of coordinates does not change the value of $\Lambda$. Let $G \in GL_2(\RR)$, and set, for all $t \geq 0$, $Z_t = G Y_t$. Then $(Z_t,I_t)$ is a PDMP with $Z$ solution of $$\frac{dZ_t}{dt}=B_{I_t} Z_t,$$ where $B_i = G A_i G^{-1}$. Due to $\lim  \frac{\log \|Y_t\|}{t}= \lim \frac{ \log \| Z_t\|}{t}$, one can see that the growth rates of $Y_t$ and $Z_t$ are equal.

Next, since the eigenvalues of $A_0$ are $\pm i \omega_0$ for $\omega_0 := \sqrt{-(\alpha_0^2 + \beta_0 \gamma_0)}$, a classical result in linear algebra (see for example \cite[Chapter 4, Theorem 3]{HS74}) states that there exists a matrix $G \in GL_2(\RR)$ such that  $$ B_0 = G A_0 G^{-1}= \begin{pmatrix}
   0 & -\omega_0 \\
   \omega_0 & 0
\end{pmatrix}.$$
Thus, replacing if necessary the matrices $A_0$ and $A_1$ by $B_0$ and $B_1 = G A_1 G^{-1}$, one may always assume that $A_0$ has the form \eqref{lem:part}.
\end{proof}

\begin{proof}[Proof of Theorem \ref{mainprop}]

We recall that it is sufficient to prove that $d=2$ and thus to show that \textbf{(i)} and \textbf{(ii)} of Theorem \ref{thm:boug} are satisfied.

We first show that \textbf{(i)} holds. For this, we need to display a matrix of $\mathcal{S}_1$ with two eigenvalues with different modulus. We claim that for every $t, s > 0$, the matrix $e^{sA_0} e^{t A_1}$ is in $\mathcal{S}_1$. Thus it is sufficient to find such a matrix with two eigenvalues with different modulus. We start by showing that it is indeed possible to find a matrix $e^{sA_0} e^{t A_1}$ for some $s, t > 0$  with two eigenvalues with different modulus, before proving the claim.

According to Lemma \ref{lem:part}, we assume  that  $A_0$ is of the form $$ A_0 = \begin{pmatrix}
   0 & -\omega_0 \\
   \omega_0 & 0
\end{pmatrix}.$$

%To prove \textbf{(i)}, we show that the matrix $$ N = \begin{pmatrix}
%   -\sqrt{\frac{\alpha_0 \beta_1}{\beta_0 \alpha_1}} & 0 \\
%   0 & -\sqrt{\frac{\alpha_1 \beta_0}{\beta_1 \alpha_0}}
%\end{pmatrix}$$
%is in $S_1$. Indeed, since by assumption $\alpha_1 \beta _0 \neq \alpha_0 \beta_1$, then $N$ has two eigenvalues with different modulus, so \textbf{(i)} will be proved.
%
%Standard computations show that for $i \in E$ and $t \in \RR$, one has $$ e^{t A_i} = \begin{pmatrix}
% \cos(\sqrt{\alpha_i \beta_i} t ) & -\sqrt{\frac{\alpha_i}{\beta_i}} \sin (\sqrt{\alpha_i \beta_i} t ) \\
%   \sqrt{\frac{\beta_i}{\alpha_i}} \sin (\sqrt{\alpha_i \beta_i} t ) & \cos(\sqrt{\alpha_i \beta_i} t )
%\end{pmatrix}.$$
%
%In particular, setting $t_i = \pi/2\sqrt{\alpha_i \beta_i}$, one has
%\beq
%N =  e^{t_0 A_0} e^{t_1 A_1}.
%\label{eq:M}
%\eeq

By standard computations, one can show that for all $t \geq 0$ and $i \in E$, $$ e^{t A_i} = \cos(\omega_i t) \mathrm{Id} + \frac{1}{\omega_i} \sin(\omega_i t) A_i,$$

where $\omega_i := \sqrt{-(\alpha_i^2 + \beta_i \gamma_i)}$. In particular, since $\mathrm{Tr}(A_i) = 0$, one has that for all $s,t \geq 0$ $$\varphi(s,t):=\mathrm{Tr}(e^{s A_0} e^{t A_1}) = 2 \cos(\omega_0 s) \cos(\omega_1 t) + \frac{1}{\omega_0 \omega_1} \sin(\omega_0s) \sin(\omega_1 t) \mathrm{Tr}(A_0 A_1).$$

On the other hand, since $\mathrm{Tr}(A_i) = 0$, one has $\det(e^{s A_0} e^{t A_1}) = 1$. Thus, denoting by $\mu_1, \mu_2$ the eigenvalues of $e^{s A_0} e^{t A_1}$ one can see that $\mu_1 \mu_2 = 1$ or equivalently $\mu_1 = 1/\mu_2$. Let assume that the eigenvalues are such that $|\mu_1| \geq |\mu_2|$. Then point \textbf{(i)} of Theorem \ref{thm:boug} is checked if $| \mu_1 | > | \mu_2 |$. This condition is equivalent to $\vert\varphi(s,t) \vert > 2$. Indeed, due to the fact that $\mu_1 + \mu_2 = \mathrm{Tr}(e^{s A_0} e^{t A_1})$, one has

$$
\vert\varphi(s,t) \vert =
\begin{cases}
\vert \mu_1 \vert + \frac{1}{\vert \mu_1 \vert} \quad \mbox{if} \ \mu_1 \in \RR\\
2\vert \mathrm{Re}(\mu_1) \vert  \quad \mbox{if} \ \mu_1 \in \CC \setminus \RR,
\end{cases}
$$
where $\mathrm{Re}(\cdot)$ stands for the real part. Combined with $\mu_1 \mu_2 = 1$, this  proves the equivalence.  By studying the derivatives of $\varphi(s,t)$, one sees that its extremal values are reached at points $(s^*, t^*)$ of the form $(n \pi/ \omega_0,m \pi/\omega_1)$ or $( \pi/2 \omega_0 + n \pi/ \omega_0 , \pi/2 \omega_1 + m \pi/\omega_1)$ with $n,m \in \ZZ$. From this we note that the extremal values are $\varphi(s^*,t^*) = \pm 2$ or $$\varphi(s^*,t^*)^2 =  \frac{1}{\omega_0^2 \omega_1^2}  \mathrm{Tr}(A_0 A_1)^2 =  \frac{( \beta_1 - \gamma_1 )^2}{\omega_1^2}, $$
where the second equality comes from the specific form of $A_0$. Therefore
\begin{equation}\label{e:equivalence}
\varphi(s^*,t^*)^2 > 4 \Longleftrightarrow \alpha_1 \neq 0~\text{or}~ \beta_1+\gamma_1 \neq 0 \Longleftrightarrow A_1 ~\text{is not proportional to}~A_0.
\end{equation}
By assumption, $A_1$ and $A_0$ are not proportional. Therefore, using \eqref{e:equivalence} one infers that the matrix $N(t_0,t_1)=N:=e^{t_0 A_0} e^{t_1 A_1}$ has two eigenvalues with different moduli for $(t_0, t_1)=(\pi/2 \omega_0,\pi/2 \omega_1)$.

In order to conclude that assumption \textbf{(i)} from Theorem \ref{thm:boug} is satisfied, we show that the matrix $ N $ lies in $\mathcal{S}_1$.

Let $V$ be a neighborhood of $N$ in $GL_2(\RR)$. Then, by  continuity, there exists $\eps > 0$ such that for all $u \in [t_0 - \eps, t_0 + \eps]$,  $s \in [t_1 - \eps, t_1 + \eps]$ and $\delta \leq \eps$, the matrix $N_{s,u,\delta}=e^{ \delta A_1} e^{u A_0} e^{s A_1}$ is in $V$. Let $V_{\eps}$ be the set of the matrices $N_{s,u,\delta}$ for $s,u$ and $\delta$ as before. Then $V_{\eps} \subset V$. Recall that $(M_t, I_t)_{t \geq 0}$ is the PDMP defined by \eqref{e:M_t}. Let $(U_n)_{n \geq 1}$ denote the sequence of interjump times of the process $I$. Then, on the event $B_{t,\eps} = \{U_1 \in [t_0 - \eps, t_0 + \eps]; U_2 \in [t_1 - \eps, t_1 + \eps]; t-(U_1+U_2) \leq \eps; U_1 + U_2 + U_3 \geq t \}$, $I_t = 1$ and $M_t \in V_{\eps}$, conditionally on $I_0 = 1$. Thus one has \begin{align*}
\PP_{\mathrm{Id}, 1} \left( (M_t,I_t) \in V \times \{1 \} \right) & \geq \PP_{\mathrm{Id}, 1} \left( (M_t,I_t) \in V_{\eps} \times \{1 \} \right)\\
& \geq \PP_{\mathrm{Id}, 1} \left( B_{t,\eps} \right).
 \end{align*}
This last probability is positive for all $\eps >0$ and $t \in [ t_0 + t_1 - 2 \eps, t_0 + t_1 +3 \eps]$. Hence $$ U((\mathrm{Id},1), V \times \{1 \}) = \int_0^{+\infty} \mathrm{e}^{-t} \PP_{\mathrm{Id}, 1} \left( (M_t,I_t) \in V \times \{1 \} \right) dt > 0.$$ This is true for all neighborhoods of $N$, so $N \in \mathcal{S}_1$ and point \textbf{(i)} is shown.

Using similar arguments, one can show that the family of matrices $\left(e^{t A_1}\right)_{t\geq 0}$ is in $\mathcal{S}_1$. Since $A_1$ has two purely imaginary eigenvalues, if $W$ denotes a finite union of one dimensional vector spaces, one has
$$
 \bigcup_{t\geq 0} e^{t A_1} W = \RR^2.
$$
In particular, for any finite union of one dimensional vector spaces $W$, there exists $t > 0$ such that $e^{t A_1} W \neq W$.  Since $e^{t A_1} \in \mathcal{S}_1$, this proves that assumption \textbf{(ii)} of Theorem \ref{thm:boug} holds.
\end{proof}

\section{Proof of Theorem \ref{mainthm} } \label{s:main}
Before we start the proof, we recall a result from \cite{BS17} that we will make use of. Let $n,m$ be two positive integers and for all $j \in W : =\{1, \ldots, n\}$, $G^j : \RR^m \to \RR^m$ be a smooth vector field such that $G^j(0) = 0$. Let $(J_t)_{t \geq 0}$ be an irreducible Markov chain on $W$ and consider the PDMP given by
$$
\frac{d \tilde X_t}{dt} = G^{J_t}(\tilde X_t).
$$
For all $j \in W$, let $B_j$ be the Jacobian matrix of $G^j$ at $0$. Like in Section \ref{sec:linear} (see also \cite{BS17} for more details), we consider the PDMP $(\Psi_t, J_t)_{t \geq 0}$ on $\mathcal{S}^{m-1} \times W$, where $\Psi$ is the angular part of the linearized process at $0$, i.e. is solution to
$$
\frac{d\Psi_t}{dt} = B_{J_t} \Psi_t - \langle B_{J_t} \Psi_t, \Psi_t \rangle \Psi_t.
$$
For $\eps >0$, let
$$
\tau^{\eps} = \inf\{ t \geq 0 \: : \| \tX_t \| \geq \eps\}.
$$Then a consequence of  \cite[Theorem 3.5, \textbf{(ii)}]{BS17} is the following :

\bthm
Assume that the PDMP $(\Psi_t, J_t)_{t \geq 0}$ admits a unique invariant probability measure $\nu$ on $\mathcal{S}^1 \times W$. If
$$
\int_{S^1 \times W} \langle B_j \psi, \psi \rangle \nu(d\psi dj) > 0,
$$
then there exist $\eps > 0$,  $\eta > 1$, $\theta > 0$ and $C > 0$ such that
for all $\bx  \in \RR^m\setminus \{0\}$ and $j \in W$, $$\mathbb{E}_{\bx,i}(\eta^{\tau^{\eps}}) \leq C (1 + \|\bx\|^{-\theta}).$$
\label{thm:lambdapositivegeneral}
\ethm

\begin{proof}
The proof of \cite[Theorem 3.5, \textbf{(ii)}]{BS17} is not given in \cite{BS17}, that is why we prove Theorem \ref{thm:lambdapositivegeneral}. In the case where there exists a compact set $K$ containing $0$ such that
$$
\tX_0 \in K \Rightarrow \tX_t \in K, \quad \forall t \geq 0,
$$
then Theorem \ref{thm:lambdapositivegeneral} is a direct consequence of \cite[Theorem 3.2, \textbf{(iii)}]{BS17}. We now show how we can still use \cite[Theorem 3.2, \textbf{(iii)}]{BS17} in the general context.

Let $K \subset \RR^d$ be a compact set containing $0$ in its interior. Let $\varphi^K : \RR^d \to [0,1]$ be a smooth function such that $\varphi^K = 1$ on $K^{\delta}$ and $\varphi^K = 0$ on the complement of $K^{2 \delta}$. Here $K^{\delta} = \{ x \in \RR^d \: : \: d(x,K) < \delta \}$ is the $\delta$ - neighbourhood of $K$.  For $i \in E$, set $G^{i,K}=\varphi^K G^i$. Note that $G^{i,K} = G^i$ on $K^{\delta}$. In particular, $0$ is a common equilibrium of the $G^{i,K}$ and $DG^{i,K}(0) = DF^i(0) = B_i$. Now consider the PDMP $(\tX^K,I)$, with $(\tX_t^K)_{ t \geq 0}$ solution of $$\frac{d \tX^K_t}{dt}=G^{I_t,K} (\tX^K_t).$$ Then we have the two following facts. First, denote by $ \tau^K = \inf \{t \geq 0 : \: \tX_t \notin K \}$ the exit time of $K$ for $\tX_t$. Then if $\tX_0 = \tX_0^K = \bx \in K$, for all $t \leq \tau^K$, $\tX_t = \tX_t^K$ almost surely. Next, since  $DG^{i,K}(0) = B_i$, the average growth rate $\Lambda^K$ of $(X^K,I)$ is equal to
$$
\Lambda^K=\int_{S^1 \times W} \langle B_j \psi, \psi \rangle \nu(d\psi dj).
$$
 Now, since $\tX_t^K$ remains in the compact set $\overline{K^{2 \delta}}$ and $\Lambda^K > 0$, one can apply \cite[Theorem 3.2, \textbf{(iii)}]{BS17}.  According to this theorem, since $\Lambda^K > 0$,  there exist $\eps > 0$, $\theta >0$,  $\eta > 1$ and $C > 0$ such that
for all $\bx \in K \setminus \{0\}$ and $i \in E$,
$$
\mathbb{E}_{\bx,i}(\eta^{\tau^{K,\eps}}) \leq C (1 + \|\bx\|^{-\theta}),
$$
 where  $\tau^{K,\eps} = \inf \{t \geq 0 : \: \|\tX_t^K\| \geq \eps \}.$ Without loss of generality, we can assume that the ball of center $0$ and radius $\eps$ is included in the interior of $K$. Let $\tau^{\eps} = \inf \{t \geq 0 : \: \|\tX_t\| \geq \eps \}.$ Now if $ \|\bx\| \geq \eps$, $\tau^{\eps} = 0$. If $\|\bx\| < \eps$, then since $\tX_t = \tX_t^K$ for all $t \leq \tau_K$, one gets that $ \tau^{\eps} = \tau^{K,\eps} \leq \tau_K$. In particular, for all $\bx \in \RR^d \setminus \{0\}$ and $i \in E$,
\begin{equation}
\mathbb{E}_{\bx,i}(\eta^{\tau^{\eps}}) \leq C (1 + \|\bx\|^{-\theta}).
\end{equation}

\end{proof}

\begin{proof}[Proof of Theorem \ref{mainthm}]
We proceed in steps.

First, we prove using Theorem \ref{mainprop} that the setting of Theorem \ref{thm:lambdapositivegeneral} applies. Then we show that it implies that $X_t$ cannot converge to $(p,q)$ and we conclude with Theorem \ref{thm:taku}.

Let $A_i$ denote the Jacobian matrix of the vector field $F^i$ at $(p,q)$, where $(p,q)$ is the common positive equilibrium of $F^0$ and $F^1$. Then $$A_i = \begin{pmatrix}
   0 & -b_i p \\
   d_i q & 0
\end{pmatrix} = \begin{pmatrix}
   0 & \beta_i \\
   \gamma_i  & 0
\end{pmatrix} ,$$
where $\beta_i = - b_i p$ and $ \gamma_i = d_i q$. The linear PDMP $(Y,I)$ where $Y$ is the solution of $$\frac{dY_t}{dt}=A_{I_t} Y_t,$$ is a particular case of the systems studied in Section \ref{sec:linear}.

To apply Theorem \ref{mainprop}, we have to check that $A_0$ and $A_1$ are non collinear. This is equivalent to showing that $\gamma_1 \beta _0 \neq \gamma_0 \beta_1$. Assume that $\gamma_1 \beta _0 = \gamma_0 \beta_1$. Then since $\beta_i = - b_i p$ and $ \gamma_i = d_i q$, we get $b_1 d_0 = b_0 d_1$. Moreover, since $p_0=p_1$ and $q_0=q_1$, one has $c_0 d_1 = c_1 d_0$ and $a_0 b_1 = a_1 b_0$.  If we set $\delta = b_1/b_0$, we note that $\kappa_1 = \delta \kappa_0$ for $\kappa \in \{ a,b, c, d \}$, which implies $F^1 = \delta F^0$. This contradicts the assumption that the vector fields $F^0$ and $F^1$ are non collinear.

As a result, $A_0$ and $A_1$ cannot be collinear. We can therefore apply Theorem \ref{mainprop} and conclude that  $\Lambda > 0$.

In particular, Theorem \ref{thm:lambdapositivegeneral} can be applied to deduce that
there exist $\eps > 0$,  $\eta > 1$, $\theta > 0$ and $C > 0$ such that
for all $\bx:=(x_0,y_0)  \in \RR_{++}^2\setminus \{(p,q)\}$ and $i \in E$,
\begin{equation}
\label{e:exit}
\mathbb{E}_{\bx,i}(\eta^{\tau^{\eps}}) \leq C (1 + \|\bx-(p,q)\|^{-\theta}).
\end{equation}

We claim that because of \eqref{e:exit} $X_t$ cannot converge to $(p,q)$. We argue by contradiction. Let $\bx \in \RR_{++}^2,  i \in E$ and assume that $X_t$ converges to $(p,q)$  almost surely under $\PP_{\bx,i}$. Define two stopping times by
$$
\tau_{\eps/2}^{\mathrm{in},1} = \inf \{ t \geq 0 \: : \: \|X_t - (p,q) \| \leq \eps / 2 \}
$$ and
$$
\tau_{\eps}^{\mathrm{out},1} = \inf \{ t > \tau_{\eps/2}^{\mathrm{in},1} \: : \: \|X_t - (p,q) \| \geq \eps  \}.
$$
Since $X_t$ converges to $(p,q)$ almost surely, one has $\PP_{\bx,i}( \tau_{\eps/2}^{\mathrm{in},1} < \infty) =1$.  Using the strong Markov property at $\tau_{\eps/2}^{\mathrm{in},1}$, one gets
$$
 \PP_{\bx,i}( \tau_{\eps}^{\mathrm{out},1} < \infty) = \EE_{\bx,i}\left( \PP_{Z_{\tau_{\eps/2}^{\mathrm{in},1}}}(\tau^{\eps} < \infty)\right)=1,
 $$
 where the second equality comes from the fact that, by \eqref{e:exit}, for all $\mathbf{y} \in \RR^2_{++} \setminus \{(p,q)\}$ and $j \in E$,
 $$
 \PP_{\mathbf{y},j}(\tau^{\eps} < \infty) =1.
 $$
  Construct recursively a family of stopping times $$\tau_{\eps/2}^{\mathrm{in},k} = \inf \{ t > \tau_{\eps}^{\mathrm{out},k-1} \: : \: \|X_t - (p,q) \| \leq \eps / 2 \}$$ and $$\tau_{\eps}^{\mathrm{out},k} = \inf \{ t > \tau_{\eps/2}^{\mathrm{in},k} \: : \: \|X_t - (p,q) \| \geq \eps  \},$$ by repeating the above procedure. Then one gets that for all $k \geq 1$, $\tau_{\eps/2}^{\mathrm{in},k}$ and $\tau_{\eps}^{\mathrm{out},k}$ are finite almost surely. This contradicts the fact that $X_t$ converges to $(p,q)$. As a result we have shown that $X_t$ cannot converge to $(p,q)$. In particular, due to Theorem \ref{thm:taku}, with probability one,
$$
\limsup x_t = \limsup y_t = +\infty, \quad \liminf x_t = \liminf y_t =0.
$$
\end{proof}
\begin{remark}\label{r:extend}
The proof of Theorem \ref{mainthm} above extends verbatim to the proof of Theorem \ref{thm:gen}.
The fact that the jump rates now depend on the position does not affect the result because when it comes to the linear system in Theorem \ref{mainprop}, one just has the constants $k_{ij}(0)$ as jump rates (see \cite[Section 2]{BS17} for details).

\end{remark}

\section{Future research}\label{s:future}

Using some of the methods developed in \cite{BS17} we were able to prove a conjecture from \cite{t06} and show that if one switches between two deterministic Lotka--Volterra systems with a common equilibrium point at $(p,q)$ then the resulting PDMP can never converge to this equilibrium.
We reduced the analysis from the non-linear Lotka--Volterra PDMP to the study of a  linear PDMP (a linearization of the original PDMP around the equilibrium point).

Recently, there have been several studies about randomly switched linear systems in dimension $2$ (see \cite{BLMZexample}, \cite{Lawley&matt&Reed} and \cite{Lag16}). In these studies, the authors show that the growth rate is positive for some switching rates by a direct computation of the invariant measure of the process $(\Theta,I)$ (this is the process that arises as the angular part when doing the polar decomposition). One could try a similar method in our setting. However, the integral expression we obtain for the growth rate does not easily yield the sign of the growth rate. Nonetheless, it could be interesting to investigate this integral expression, possibly through numerical simulations.

Another interesting direction for the future is finding out whether the process $X_t$ defined by \eqref{e:LV_PDMP} is transient, null-recurrent of positive recurrent. The simulations done in \cite{t06} seem to suggest the following conjecture.

\begin{conjecture}
Suppose $X_t=(x_t,y_t)$ is the process defined by \eqref{e:LV_PDMP} together with the initial condition $X_0=(x_0,y_0)\in\RR_{++}^2$. Then, almost surely
\[
\lim_{t\to\infty} \left( x_t+y_t+\frac{1}{x_t}+\frac{1}{y_t}\right)=\infty
\]
and the process $X_t$ is transient.
\end{conjecture}

{\bf Acknowledgments.} The authors thank Michel Bena{\"\i}m and Dang Nguyen for very helpful discussions. The comments of two anonymous reviewers have significantly improved the content and presentation of the paper.
\bibliographystyle{siam}
\bibliography{RandomSwitch}
\end{document}